\documentclass[12pt,epsf]{article}
\usepackage{amsmath,amssymb}
\usepackage[dvips]{graphicx}

\newtheorem{thm}{Theorem}[section]
\newtheorem{prop}[thm]{Proposition}
\newtheorem{lem}[thm]{Lemma}
\newtheorem{df}[thm]{Definition}

\begin{document}
\title{The automorphism group of the singular\\ $K3$ surface of discriminant $7$}
\author{Masashi Ujikawa}
\maketitle

\begin{abstract}
We give a system of generators of the automorphism group of the singular $K3$ surface of discriminant $7$. 
This system of generators consists of the inversion involutions of some elliptic fibrations with a section together with ${\rm PGL}_2(7)$.
\end{abstract}


\section{Introduction}

In this paper, we work over the field of complex numbers.
We restrict the term {\it $K3$ surface} to the algebraic one, i.e. a $K3$ surface is a simply-connected {\it algebraic} surface with the trivial canonical line bundle.  

The automorphism group of a $K3$ surface has been studied for a long time. 
The Torelli theorem for $K3$ surfaces by Piatetski-Shapiro and Shafarevich \cite{PS-S} states that the automorphism group of a $K3$ surface is isomorphic to the group of isomorphisms of the second cohomology group which preserve the intersection product and the ample cone and fix the ray of holomorphic $2$-forms. 
Shioda and Inose \cite{S-I} obtained the result that any singular $K3$ surface has an infinite group of automorphisms. 

A $K3$ surface is said to be singular if its Picard number is $20$ (maximal possible). 
It is interesting to calculate a system of generators of the automorphism group of a singular $K3$ surface.
Vinberg \cite{2most} calculated the automorphism groups of two singular $K3$ surfaces. 
These singular $K3$ surfaces are of discriminants $3$ and $4$. 
The purpose of this paper is to give a system of generators of the automorphism group of the singular $K3$ surface with discriminant $7$, the next small discriminant. 
For the cases of larger discriminants, Keum and Kond{\=o} \cite{KK} gave systems of generators of the automorphism groups of the singular $K3$ surfaces with discriminants $12$ and $16$. 
The transcendental lattices of these singular $K3$ surfaces with discriminants $3$, $4$, $7$, $12$, $16$ have Gram matrices 
\[
\left[\begin{array}{cc}2&1\\1&2\end{array}\right],
\left[\begin{array}{cc}2&0\\0&2\end{array}\right],
\left[\begin{array}{cc}2&1\\1&4\end{array}\right],
\left[\begin{array}{cc}4&2\\2&4\end{array}\right],
\left[\begin{array}{cc}4&0\\0&4\end{array}\right],
\]
respectively. 

Here we explain the method of calculating a system of generators of the automorphism group in this paper. 
This method was used for the calculations of the automorphism groups of $K3$ surfaces in several papers \cite{Kondo}, \cite{KK}, \cite{Hessian}, \cite{p2s1}, \cite{p3s1}.

The automorphism group of a $K3$ surface $Y$ is the factor group of the orthogonal group of the N{\'e}ron-Severi lattice $S_Y$ of $Y$ by the reflection group generated by the $(-2)$-roots of $S_Y$ up to finite groups, and the difference can be described explicitly (Theorem \ref{Torelli}). 
The ample cone $D(S_Y)$ of $Y$ is a fundamental domain of the $(-2)$-reflection group, and the above factor group is the automorphism group ${\rm Aut}(D(S_Y))$ of the ample cone. 
In case of the singular $K3$ surface $X$ of discriminant $7$, it is difficult to describe ${\rm Aut}(D(S_X))$ directly, since $D(S_X)$ has infinitely many faces. 

Instead, we embed the N{\'e}ron-Severi lattice $S_X$ to the even unimodular lattice $I{\!}I_{1,25}$ of signature $(1,25)$ as the orthogonal complement of a root lattice $A_6$. 
The isometries of $S_X$ extend to those of the ambient lattice $I{\!}I_{1,25}$. 
Therefore the orthogonal group of $S_X$ is the subgroup of the orthogonal group of $I{\!}I_{1,25}$ preserving the root lattice $A_6$. 
On the other hand, Conway \cite{Con-Lor} found a fundamental domain $D$ of the reflection group of $I{\!}I_{1,25}$ whose faces bijectively correspond to the Leech lattice points (Theorem \ref{Conway}). 
Borcherds used Conway's fundamental domain to find explicit descriptions of the automorphism groups of some Lorentzian lattices \cite{Bor-Lor}. 
Borcherds also treated the automorphism group of the N{\'e}ron-Severi lattice of the singular $K3$ surface $X$ of discriminant $7$ \cite[Example 5.4.]{Bor-K3}. 
In this situation, the restriction $D'$ of the fundamental domain $D$ to $S_X{\otimes}{\mathbb R}$ is contained in the ample cone $D(S_X)$. 
Furthermore, the polyhedral cone $D'$ has only a finite number of faces. 
The automorphism group of $D^{\prime}$ is isomorphic to ${\rm PGL}_2(7)$. 
The faces of $D^{\prime}$ consist of four types, namely $A_6\oplus A_1$, $A_7$, $D_7$ and $E_7$ (Proposition \ref{add-root}, Table \ref{face-tab}). 
The Leech roots giving the faces of type $A_6\oplus A_1$ correspond to twenty-eight smooth rational curves on $X$. 
For each of the faces of the other types, there exists an elliptic fibration whose fiber class is contained in the face. 
The inversion involutions associated to these elliptic fibrations play an analogous role of the reflections with respect to the corresponding faces (see Proposition \ref{q-refl} for details). 
It can be proved that these inversion involutions together with ${\rm PGL}_2(7)$ generate the automorphism group of the ample cone (Theorem \ref{gen}), and the procedure is completed. 

This paper is organized as follows. 
In Section 2, we recall some known results on lattices, $K3$ surfaces and elliptic fibrations. We also introduce Conway's fundamental domain. 
In Section 3, we recall some results obtained by Naruki \cite{Naruki}.
In Section 4, we identify the singular $K3$ surface $X$ of discriminant $7$ with Naruki's $K3$ surface, and we show that the automorphism group of $X$ is isomorphic to the automorphism group of the ample cone $D(S_X)$ of $X$. 
Then we embed the N{\'e}ron-Severi lattice of the singular $K3$ surface of discriminant $7$ into $I{\!}I_{1,25}$, and  study the restriction of Conway's fundamental domain to the N{\'e}ron-Severi lattice. 
In Section 5, we give some elliptic fibrations of the singular $K3$ surface $X$ of discriminant $7$ and consider their inversion involutions. 
These involutions together with ${\rm PGL}_2(7)$ are sufficient to generate the automorphism group of the singular $K3$ surface $X$ of discriminant $7$, as will be proved in Section 6.

\vspace{2mm}
{\it Acknowledgements.} \quad The author expresses his gratitude to Professor Shigeyuki Kond{\=o} for suggesting this problem, helpful discussions and comments. 
He is also grateful to Professor Hisanori Ohashi for valuable advices. 
Finally, he wishes to thank the referee for his numerous helpful comments. 
This research is partially supported by JSPS-grant 22224001.


\section{Lattices, $K3$ surfaces and elliptic fibrations}

A {\it lattice} $L$ is a free $\mathbb{Z}$-module of finite rank with a non-degenerate symmetric bilinear form 
\[
\langle \; ,\rangle : L\times L \to \mathbb{Z}.
\]
The coefficients extension $L\otimes_{\mathbb Z} R$ is denoted by $L_R$, e.g. $L_{\mathbb Q}$ or $L_{\mathbb R}$.

A lattice $L$ is said to be {\it even} if $\langle x,x\rangle$ is even for every element $x$ of $L$, and the {\it discriminant} $d(L)$ of $L$ is the determinant of the Gram matrix ${\rm Gram}_L$ of a basis of $L$. 
A lattice $L$ of rank $\rho$ is said to be {\it Lorentzian} if its signature is $(1,{\rho}-1)$. 
The orthogonal group of $L$ is denoted by ${\rm O}(L)$.

For a lattice $L$, ${\rm Hom}_{\mathbb Z}(L,{\mathbb Z})$ is called the {\it dual} lattice of $L$. 
The dual lattice of $L$ is denoted by $L^{\vee}$. 
The dual lattice $L^{\vee}$ can be regarded as a submodule of $L_{\mathbb Q}$. 
A {\it root} of $L$ is an element $r$ of $L$ such that $2r/\langle r,r\rangle$ is contained in $L^{\vee}$. 
A root $r$ is called an $n$-root if $r^2=n$.
For a root $r$, the {\it reflection} $s_r\in {\rm O}(L)$ is defined to be 
\[
s_r : x\mapsto x-2\frac{\langle x,r\rangle }{\langle r,r\rangle }r.
\]
Even if $r$ is a vector in $L_{\mathbb Q}$ which is not a root, $s_r$ denote the ${\mathbb Q}$-linear transformation on $L_{\mathbb Q}$ defined by the same rule. 

\vspace{2mm}
Let $L$ be a Lorentzian lattice. Then $\{ x\in L_{\mathbb R}\mbox{ : }x^2>0\}$ consists of two connected components. We denote one of them by $P$. 
We denote by $W(L)$ (resp. $W(L)^{(2)}$) the subgroup of ${\rm O}(L)$ generated by the reflections with respect to the roots (resp. $(-2)$-roots) $r$ such that $r^2<0$. 
These are normal subgroups of ${\rm O}(L)$, because $gs_rg^{-1}=s_{g(r)}$ for any element $g$ of ${\rm O}(L)$. 
A connected component ${\rm C}$ of the complement $P\setminus \bigcup r^{\perp}$ of all reflection hyperplanes is called a {\it Weyl chamber} whose closure $\overline{\rm C}$ is a fundamental domain of the action of $W(L)$ on $P$. 
If a Weyl chamber ${\rm C}$ is given, a root $r$ is said to be {\it simple} if $r$ corresponds to a hyperplane of the boundary of ${\rm C}$ such that $\langle r,c\rangle >0$ for a point $c$ in ${\rm C}$.

\vspace{2mm}
A {\it root lattice} is a {\it negative}-definite lattice generated by $(-2)$-roots. 
This is a direct sum of the irreducible root lattices which are classified into the types $A_l$, $D_m$, $E_n$, where $l\geq 1$, $m\geq 4$, $n=6$, $7$ and $8$. 

We denote by $U$ the {\it hyperbolic plane}, that is, the lattice of rank $2$ with the Gram matrix 
\[
\left[ \begin{array}{cc}0&1\\1&0\end{array}\right] .
\] 
This is the unique even unimodular lattice of signature $(1,1)$. 

We denote by $\Lambda$ the {\it Leech lattice}, i.e. the unique even unimodular lattice of signature $(0,24)$ without any $(-2)$-roots \cite{CS}. 
Here, we describe a construction of the Leech lattice. 
First, let $\Omega$ be the projective line over the finite field ${\mathbb F}_{23}$:
\[
\Omega =\{ \infty , 0,1,\ldots , 22\} .
\]
A Steiner system $S(5,8,24)$ is a set of $8$-elements subsets of $\Omega$ such that, for arbitrary $5$-elements subset $\alpha$ of $\Omega$, there exists a unique element $\beta$ of $S(5,8,24)$ containing $\alpha$. 
An element of a Steiner system is called an {\it octad}. 
We use the Steiner system constructed by Carmichael \cite{Carmichael}, which is the ${\rm PSL}_2(23)$-orbit of $\{$ $\infty$ $0$ $1$ $3$ $12$ $15$ $21$ $22$ $\}$ (listed by Todd \cite{Todd}). (This is the Steiner system used by Conway and Sloane \cite{CS}). 
In the sequel, we use this Steiner system.

A Steiner system $S(5,8,24)$ generates in ${\mathbb F}_2^{\Omega}$ a $12$-dimensional linear subspace ${\cal C}$. 
This subspace is called the {\it (binary) Golay code}. 
We call a supporting subset of an element of ${\cal C}$ a ${\cal C}$-set.
\begin{df}{\rm 
The Leech lattice $\Lambda$ is the set of all integral vector points $(\xi_{\infty}$, $\xi_0$, $\ldots$, $\xi_{22})$ in ${\mathbb Z}^{24}$ which satisfy the following conditions:  
\begin{enumerate}
\item There exists $m$, such that the coordinates $\xi_i$ are congruent to $m$ modulo $2$ for any $i$.
\item For any value $n$, the set of $i$ for which $\xi_i$ is congruent to $n$ modulo $4$ is a ${\cal C}$-set.
\item The coordinate sum $\sum \xi_i$ is congruent to $4m$ modulo $8$.
\end{enumerate}
The bilinear form on $\Lambda$ is defined for two elements $\xi =(\xi_{\infty},\xi_0,\ldots ,\xi_{22})$, $\eta =(\eta_{\infty},\eta_0,\ldots ,\eta_{22})$ in $\Lambda$, to be 
\[
\langle {\xi},{\eta}\rangle{=}-\frac{\sum \xi_i\eta_i}{8}.
\]
}
\end{df}
For a subset $\alpha$ of $\Omega$, we denote by $\nu_{\alpha}$ the vector in ${\mathbb R}^{24}$ with entries $1$ on $\alpha$ and entries $0$ on $\Omega\setminus\alpha$.

Let $I{\!}I_{1,25}$ be the unique even unimodular lattice of signature $(1,25)$. Then it is isometric to $U\oplus \Lambda$. 
With respect to this decomposition, hereafter we denote a lattice vector in $I{\!}I_{1,25}$ by $(m,n;\lambda )$ where $(m,n)$ is the element in $U$ and $\lambda$ is the element in $\Lambda$. 
If a lattice vector in $I{\!}I_{1,25}$ is of the form $(-\lambda^2/2-1,1;\lambda )$ for some $\lambda\in\Lambda$, then it is a $(-2)$-root, and this is called a {\it Leech root}. 

If a Weyl chamber ${\rm C}$ is given, a lattice vector $w$ satisfying $\langle w,r_i\rangle = -r_i^2/2$ for any simple root $r_i$ is called a {\it Weyl vector}. 
The vector $(1,0;0) \in U\oplus\Lambda$ is a Weyl vector of Conway's fundamental domain by the following. 

\begin{thm}{\rm (Conway \cite{Con-Lor})}\label{Conway}
Define a polyhedral cone $D$ in the positive cone $P{\subset}(I{\!}I_{1,25})_{\mathbb R}$ by
\[
D=\{ x\in P\, :\,\,\langle x,r\rangle >0 \mbox{ for any Leech root }r\} .
\]
Then $D$ is a fundamental domain of $W(I{\!}I_{1,25})^{(2)} = W(I{\!}I_{1,25})$. 
Furthermore, the orthogonal group ${\rm O}(I{\!}I_{1,25})$ is a split extension of ${\rm Aut}(D)$ by $\{\pm{1}\} {\cdot}W(I{\!}I_{1,25})$. 
Here ${\rm Aut}(D)$ is the stabilizer of $D$ in ${\rm O}(I{\!}I_{1,25})$, which is isomorphic to the affine automorphism group $\cdot\infty$ of $\Lambda$.
\end{thm}

Let $L$ be an even lattice. The {\it discriminant group} $A_L$ of $L$ is defined to be the abelian group $L^{\vee}/L$ of order $|d(L)|$.  
Then, the quadratic form $x\mapsto\langle x,x\rangle$ on $L$ gives  a finite quadratic form $q_{L}$ : $ A_L\to \mathbb{Q}/2\mathbb{Z}$ on $A_L$, i.e. $q_L$ satisfies
\begin{enumerate}
\item $q_{L}(na)=n^2q_{L}(a)$ for any $n\in \mathbb{Z}$, $a\in A_L$,
\item $q_L(a+a')-q_L(a)-q_L(a')$ is the duplication of a symmetric bilinear form on $A_L$ with value in $\mathbb{Q}/\mathbb{Z}$.
\end{enumerate}
The finite quadratic form $q_L$ is called the {\it discriminant quadratic form}.
The automorphism group of $A_L$ preseving $q_L$ is denoted by ${\rm O}(q_L)$. 
There exists a natural map ${\rm O}(L)$ $\to$ ${\rm O}(q_L)$. 

\begin{thm}{\rm (Nikulin \cite[Proposition 1.4.2]{Nik-lat})}\label{Nikulin}
Let $L$ be an even unimodular lattice, let $S$ be a primitive sublattice of $L$ and let $T$ be the orthogonal complement of $S$ in $L$. 
Then $A_S$ and $A_T$ is anti-isometrically isomorphic. 

Furthermore, let $g$ be an element of ${\rm O}(S)$ and let $h$ be an element of ${\rm O}(T)$. 
Then $g{\oplus}h$ is an element of ${\rm O}(L)$ if and only if the image of $g$ by ${\rm O}(S)$ $\to$ ${\rm O}(A_S)$ and the image of $h$ by ${\rm O}(T)$ $\to$ ${\rm O}(A_T)$ coincide. 
\end{thm}

For a $K3$ surface $Y$, the second cohomology group $H^2(Y,{\mathbb Z})$ with the cup product is an even unimodular lattice of signature $(3,19)$. 
Hence $H^2(Y,{\mathbb Z})$ is isomorphic to $U^{\oplus 3}\oplus E_8^{\oplus 2}$. 
For a $K3$ surface $Y$ with a nowhere vanishing holomorphic $2$-form $\omega_Y$, the {\it N{\'e}ron-Severi lattice} $S_Y$ of $Y$ is defined by 
\[
S_Y=\{ x\in H^2(Y,{\mathbb Z})\mbox{ : }\langle x,\omega_Y\rangle =0\} .
\]
The orthogonal complement $T_Y$ of $S_Y$ in $H^2(Y,{\mathbb Z})$ is called the {\it transcendental lattice} of $Y$. 
Since $H^2(Y,{\mathbb Z})$ is unimodular, the discriminant groups of $S_Y$ and $T_Y$ are anti-isometrically isomorphic. 

Let $P(S_Y)$ be the connected component of the cone $\{ x\in (S_{Y})_{\mathbb R}$ : $\langle x$, $x\rangle >0\}$ containing an ample class. 
Then the ample cone $D(S_Y)$ of $Y$ is defined by 
\[
\begin{split}
D(S_Y)=\{ x\in P(S_Y)\mbox{ : }&\langle x,\delta \rangle >0\\ &\mbox{ for any class $\delta$ of smooth rational curves on }Y\} .
\end{split}
\]
This polyhedral cone is a fundamental domain of $W(S_Y)^{(2)}$. 

\begin{thm}\label{Torelli}{\rm (The Torelli theorem for algebraic $K3$ surfaces \cite{PS-S})}
Let $Y$, $Y^{\prime}$ be two algebraic $K3$ surfaces, and let $\varphi$ be an isometry $H^2(Y^{\prime},\mathbb{Z})$${\to}$$H^2(Y,\mathbb{Z})$ which satisfies the conditions (i) ${\varphi}(D(S_Y^{\prime}))$ $\subset$ ${\varphi}(D(S_Y))$ and (ii) ${\varphi}_{\mathbb C}({\omega}_{Y^{\prime}})$ $\in$ ${\mathbb C}{\omega}_Y$, where ${\varphi}_{\mathbb C}$ $={\varphi}{\otimes}{\mathbb C}$.

Then $\varphi$ is induced by an isomorphism $Y\to Y^{\prime}$.
\end{thm}

By the Torelli theorem of $K3$ surfaces, ${\rm Aut}(Y)$ is isomorphic to 
\[
{\rm O}(S_Y)/\{\pm{1}\}{\cdot}W(S_Y)^{(2)} \, {\cong}\, {\rm Aut}(D(S_Y))
\]
up to finite groups. 
More precisely, under the anti-isomorphism $A_{S_Y}\cong A_{T_Y}$, ${\rm Aut}(X)$ is isomorphic to the fiber product 
\[
{\rm Aut}(D(S_Y))\underset{{\rm O}(q_{S_Y})}{\times}H_Y
\]
where $H_Y$ is the stabilizer of ${\mathbb C}\omega_Y$ in ${\rm O}(T_Y)$, which is finite.
Now, we recall the facts \cite[Theorem 3.1]{Nik-gp} that ${\rm Aut}(Y)$ acts on $H^0(Y,\Omega^2_Y )$ $\cong$ ${\mathbb C}\omega_Y$ as a finite cyclic group, and that ${\rm O}(T_Y)$ acts on ${\mathbb C}\omega_Y$ effectively. 
Hence $H_Y$ is a finite cyclic group. 
In particular, ${\rm Aut}(Y)$ is an extension of the cyclic group ${\rm ker}$$(H_Y$$\to$${\rm O}(q_{S_Y}))$ by a finite-index subgroup of ${\rm Aut}(D(S_Y))$. 

\vspace{2mm}
We recall general notions on elliptic fibrations. 
In this paper, an {\it elliptic fibration} means an elliptic fibration with a section, i.e. a triple $(Y, f, O)$ where $Y$ is a surface, $f$ is a morphism from $Y$ to a curve $C$ whose general fiber is an elliptic curve and $O$ is a section of $f$. 
We call $O$ the {\it zero section} of this elliptic fibration. 
The zero element of each general fiber is the intersection of the fiber with the zero section. 
There exists a rational map sending a point in a general fiber to its inverse element. 
In the case that $Y$ is a $K3$ surface, this rational map extends to the whole surface by the minimality of $K3$ surfaces. 
We call this automorphism an {\it inversion involution}. 

A fiber of an elliptic fibration not isomorphic to an elliptic curve is called a {\it singular fiber}. 
Types of singular fibers are classified by Kodaira \cite[Theorem 6.2]{Kod}. 
For each singular fiber $F_i$, the irreducible components $\Theta_{i,j}$ of $F_i$ disjoint from the zero section $O$ generate a root lattice $\Theta_i=\oplus_j{\mathbb Z}\Theta_{i,j}$ in $S_Y$. 
A singular fiber has a group structure if some appropriate points are removed \cite[Theorem 9.1]{Kod}.
When this root lattice is of type $A_{n-1}$, the corresponding singular fiber has the group structure ${\mathbb C}^{\times}\times ({\mathbb Z}/{n\mathbb Z})$.
A singular fiber of type $D_n$ has the group structure ${\mathbb C}\times ({\mathbb Z}/{2\mathbb Z})^2$ or ${\mathbb C}\times ({\mathbb Z}/{4\mathbb Z})$ according to whether $n$ is even or odd.
A singular fiber of type $E_n$ has the group structure ${\mathbb C}\times ({\mathbb Z}/{(9-n)\mathbb Z})$.

Hence, the inversion involution acts on the root lattice $\Theta_i$ as the opposition involution of the root lattice (defined as the composition of the $(-1)$-multiplication and the element of the maximum word-length in the corresponding Coxeter group). 
The opposition involution is trivial for $A_1$, $D_n$ (if $n$ is even), $E_7$, and $E_8$. 
The opposition involution of $A_n$ (if $n>1$) is\\
\begin{center}
\includegraphics[height=30pt]{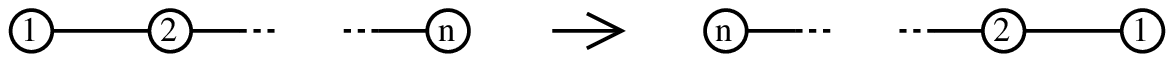}.\\
\end{center}
The opposition involution of $D_n$ (if $n$ is odd) is 
\begin{center}
\includegraphics[height=60pt]{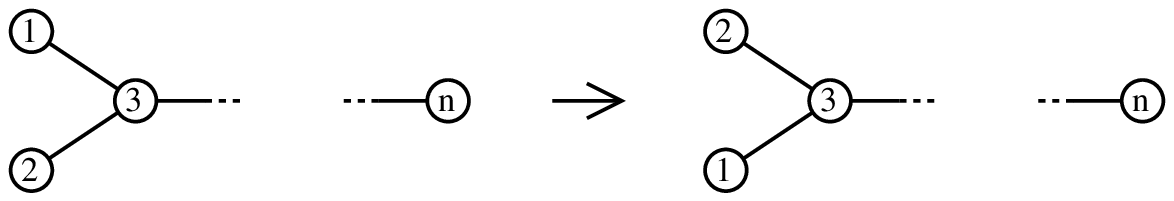}.\\
\end{center}
The opposition involution of $E_6$ is
\begin{center}
\includegraphics[height=60pt]{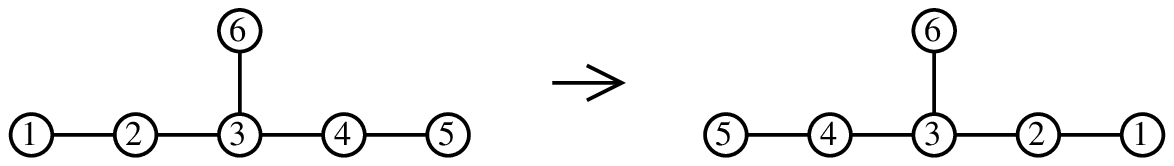}.\\
\end{center}

In the N{\'e}ron-Severi lattice of a surface with an elliptic fibration $f$, the class of a general fiber $F$ and the zero-section $O$ span a hyperbolic plane. We denote this by $U_f$:
\[
U_f=\langle F,O\rangle .
\]
The {\it trivial lattice} ${\rm Triv}$ of an elliptic fibration is the sublattice of the N{\'e}ron-Severi lattice generated by the hyperbolic plane $U_f$ and the classes of the irreducible components $\Theta_{i,j}$ of the singular fibers:
\[
{\rm Triv}=U_f\oplus \left(\bigoplus_i\Theta_i\right) .
\]

The sections of an elliptic fibration $f$ form a finitely generated abelian group. 
This group is called the {\it Mordell-Weil group}, and denoted by $MW(f)$. 
Then the following isomorphism holds \cite{MW}:
\[
MW(f)\cong S_Y/{\rm Triv}.
\]
The orthogonal complement of ${\rm Triv}$ in $S_Y$ is called the {\it essential lattice} of the elliptic fibration, and is denoted by ${\rm Ess}$. 
In particular, the orthogonal sum ${\rm Triv}\oplus{\rm Ess}$ is of finite index in $S_Y$.

\begin{prop}\label{inv}
The inversion involution of an elliptic fibration $(Y$, $f$, $O)$ acts as the $(-1)$-multiplication on ${\rm Ess}$. 
The inversion involution acts on the hyperbolic plane $U_f$ as the identity, and it acts as the opposition on the root lattice spanned by components of a singular fiber. 
\end{prop}
(Proof.)
The inversion involution acts as an isometry on N{\'e}ron-Severi group $S_Y$.
On the other hand, it acts as $-1$ on $MW(f)\cong S_Y/{\rm Triv}$.
Therefore it acts as $-1$ on ${\rm Ess}$.
\hfill $\Box$ 


\section{Naruki's results}\label{Naruki-results}

Here we give a sketch of the results of Naruki \cite{Naruki}. 

We denote ${\rm exp}(2{\pi}\sqrt{-1}/7)$ by $\zeta$, the integer ring of the cyclotomic field $\mathbb{Q}(\zeta )$ by $\mathfrak{o}$, and let $\mathfrak{p}$ be the maximal ideal of $\mathfrak{o}$ generated by $1-\zeta$.

Naruki defined a Hermitian metric $H_7$ of signature $(2,1)$ on ${\mathbb C}^3$: 
\[
H_7(z_1,z_2,z_3)=z_1\overline{z_1}+z_2\overline{z_2}-({\zeta}+\overline{\zeta})z_3\overline{z_3}.
\]
Then, the domain in ${\mathbb P}^2({\mathbb C})$
\[
B=\{ (z_1:z_2:z_3)\mbox{ : }H_7(z)<0\}
\]
is biholomorphically equivalent to the complex $2$-ball, on which $SU(2,1)$ acts.

Let $\Gamma_7$ be the discrete subgroup $SU(H_7)(\mathfrak{o})$ of special unitary group $SU(H_7)$ and $\Gamma_7'$, $\Gamma_7''$ the congruent subgroups of $\Gamma_7$ modulo $\mathfrak{p}$, $\mathfrak{p}^2$.

Naruki's main result is as follows: 

\begin{thm}{\rm (Naruki \cite{Naruki})}
The quotient $B/\Gamma_7'$ is a smooth $K3$ surface, which is singular of discriminant $7$. 
The branch locus of $B/\Gamma_7''$ $\to$ $B/\Gamma_7'$ is twenty-eight smooth rational curves, and the group $\Gamma_7/\Gamma_7'$ $\cong$ ${\rm PSL}_2(7)$ acts as permutations of these curves. 
\end{thm}


\begin{figure}[t]
\caption{Coxeter's graph. (The vertices are numbered according to the corresponding octads in the text.)}
\begin{center}
\includegraphics[height=240pt]{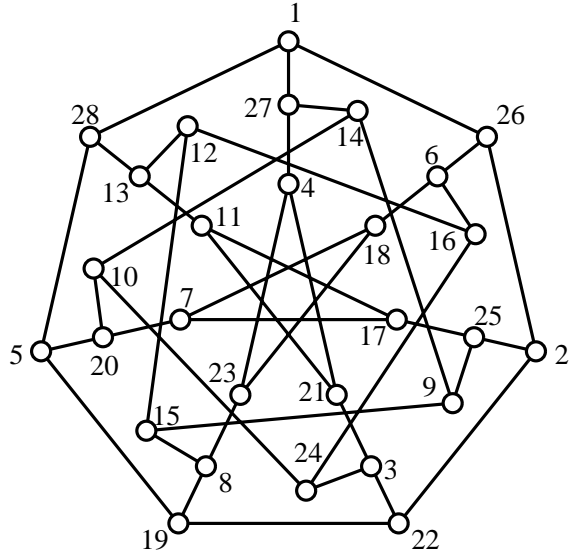}\\
\end{center}
\label{Cox-graph}
\end{figure}


Furthermore, Naruki decomposed the map $B/\Gamma_7''\to B/\Gamma_7'$ as 
\[
B/\Gamma_7''{\to}S(7)^{\sim}{\to}S(7){\to}B/\Gamma_7'
\]
where $S(7)$ is the elliptic modular surface \cite{EMS} of level $7$, and $S(7)^{\sim}$ is a cyclic cover of $S(7)$ constructed by Livn{\'e} \cite{Livne} and Inoue \cite{Inoue}. 
This surface $S(7)^{\sim}$ is of general type. 

The elliptic modular surface $S(7)$ has an elliptic fibration $S(7){\to}X(7)$, where $X(7)$ is the elliptic modular curve of level $7$, which is isomorphic to the Klein quartic curve 
\[
{\{}X^3Y+Y^3Z+Z^3X=0{\}}\subset{\mathbb P}^2 ,
\]
of genus $3$. 
This elliptic fibration has $49$ sections and the Mordell-Weil group is isomorphic to $({\mathbb Z}/7{\mathbb Z})^2$. 
The elliptic fibration has the twenty-four singular fibers all of which are of type $I_7$. 
The automorphism group of $S(7)$ is $({\mathbb Z}/7{\mathbb Z})^2{\!}:{\!}{\rm SL}_2(7)$, where $({\mathbb Z}/7{\mathbb Z})^2$ is the translations by the $49$ torsion sections, and the colon $:$ denotes a semidirect product. 

Let $N(7)$ be the normalizer of $\pi_1(S(7)^{\sim})$ in ${\rm SU}(H_7)$. Then  $N(7)$ is a extension of ${\rm Aut}(S(7)^{\sim})$ by $\pi_1(S(7)^{\sim})$. 
Livn{\'e} \cite{Livne} gave a system of generators and relations of $N(7)$ explicitly. 
Upon Livn{\'e}'s result, Naruki characterized $N(7)$ as the automorphism group of an $\mathfrak{o}$-lattice as follows. 
Let $L$ be the $\mathfrak{o}$-module $\mathfrak{o}e_1+\mathfrak{o}e_2+\mathfrak{o}e_3$ of rank $3$. Then $\Gamma_7$ is 
\[
\{ M\in {\rm SU}(H_7)\mbox{ : }ML\subseteq L\} .
\]
Let $\tilde{e_2}$ be $(\zeta -\zeta^{-1})^{-1}e_2$ $+(\zeta^3-\zeta^{-3})^{-1}e_3$, and let $L^{\sim}$ be the $\mathfrak{o}$-module $\mathfrak{o}e_1+\mathfrak{o}\tilde{e_2}+\mathfrak{o}e_3$. Then $\mathfrak{p}L^{\sim}\subseteq L\subseteq L^{\sim}$ hold. 
\begin{thm}
$N(7)$ $=$ $\{ M\in {\rm SU}(H_7)$ : $M L^{\sim}$ $\subseteq$ $L^{\sim}\}$.
\end{thm}

From these considerations, the modulo $\mathfrak{p}$ reduction of $N(7)$ is isomorphic to the finite affine transformation group $({\mathbb Z}/7{\mathbb Z})^2{\!}:{\!}{\rm SL}_2(7)$ of the sections of $S(7)$, and the reduction coincides with the projection of $N(7)$ to ${\rm Aut}(S(7))$. 
Furthermore, the modulo $\mathfrak{p}$ reduction of $\Gamma_7^{\prime}$ is the subgroup of $({\mathbb Z}/7{\mathbb Z})^2{\!}:{\!}{\rm SL}_2(7)$ isomorphic to $({\mathbb Z}/7{\mathbb Z}){\!}:{\!}({\mathbb Z}/7{\mathbb Z})$. 
Hence $B/\Gamma_7^{\prime}$ is the quotient of $S(7)$ by $({\mathbb Z}/7{\mathbb Z}){\!}:{\!}({\mathbb Z}/7{\mathbb Z})$. 


\begin{figure}[t]
\caption{A rough sketch of the real locus of the sextic curve, which is the branch locus of the morphism of degree $2$ from $B/\Gamma_7^{\prime}$ to the plane.}
\begin{center}
\includegraphics[height=240pt]{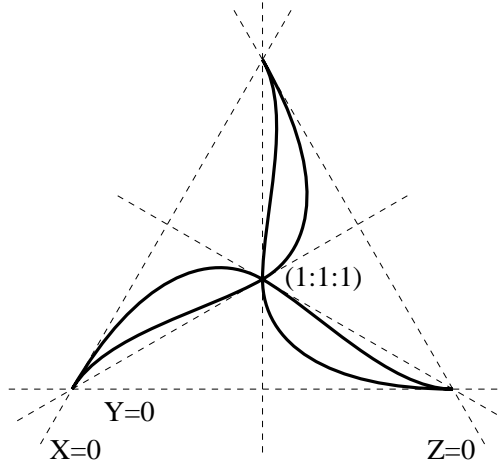}\\
\end{center}
\label{sextic}
\end{figure}


The normal subgroup ${\mathbb Z}/7{\mathbb Z}$ of $({\mathbb Z}/7{\mathbb Z}){\!}:{\!}({\mathbb Z}/7{\mathbb Z})$ acts on $S(7)$ as translations, fixes the three singular fibers, and permute the irreducible components of each of the remaining twenty-one singular fibers. 
The quotient singularities disappear by the quotient formation by the other ${\mathbb Z}/7{\mathbb Z}$. 
The action of $({\mathbb Z}/7{\mathbb Z}){\!}:{\!}({\mathbb Z}/7{\mathbb Z})$ induces the action of ${\mathbb Z}/7{\mathbb Z}$ on $X(7)$ with three fixed points. 
Hence the quotient $X(7)/({\mathbb Z}/7{\mathbb Z})$ is isomorphic to ${\mathbb P}^1$. 
Thus we obtained the commutative diagram
\begin{eqnarray*}
S(7)&\to&X(7)\\
\downarrow&&\downarrow\\
B/\Gamma_7'&\to&{\mathbb P}^1
\end{eqnarray*}
where the vertical arrows are the quotient maps, and the row arrows are the elliptic fibrations. 
The bottom row is the elliptic fibration $f_1$ with the singular fibers of type $3I_7+3I_1$. 
Twenty-one of the twenty-eight branch curves are the irreducible components of the singular fibers of type $I_7$, and the remaining seven are the sections. 

Naruki described the configuration of twenty-eight curves in the argument similar to that of Coxeter \cite{Cox} (Figure \ref{Cox-graph}). 

Blowing down $(-1)$-curves of the quotient of $B/\Gamma_7^{\prime}$ by the inversion involution $\iota_1$ of the fibration $f_1$ (similar to the argument in \cite{Naruki-pre}), Naruki obtained the double plane model of $X$ branching along the rational sextic curve with ($D_4+3A_4$)-singularities
\[
(X^2Y+Y^2Z+Z^2X-3XYZ)^2-4XYZ(X-Y)(Y-Z)(Z-X)=0
\]
drawn in Figure \ref{sextic}. 

Recently Harrache and Lecacheux \cite{HL} studied elliptic fibrations on this surface.


\section{The N{\'e}ron-Severi lattice of the singular $K3$ surface of discriminant $7$}

In this Section, we first describe the N{\'e}ron-Severi lattice $S_X$ of the singular $K3$ surface of discriminant $7$.
Then we study the restriction of Conway's fundamental domain to $(S_X)_{\mathbb R}$. 
Especially, we describe the four types of its faces.

For a singular $K3$ surface $Y$, $T_Y$ is a positive-definite lattice of rank two. 
A basis $e_1$, $e_2$ of $T_Y$ is said to be oriented if the imaginary part of $\langle e_1$, $\omega_Y\rangle$ $/$ $\langle e_2$, $\omega_Y\rangle$ is positive. 
We denote the Gram matrix of $T_Y$ with respect to a oriented basis by ${\rm Gram}_{T_X}$. 

\begin{thm}{\rm (Shioda-Inose \cite{S-I})}
There is a natural bijective correspondence from the set of singular $K3$ surfaces to the set of equivalence classes of positive-definite even integral binary quadratic forms with respect to $SL_2({\mathbb Z})$
\[
Y{\mapsto}\left[ {\rm Gram}_{T_Y}\right] .
\]
\end{thm}

Since positive-definite even binary quadratic forms of discriminant $7$ constitute only one class i.e. the class of 
\[
\left[\begin{array}{cc}2&1\\1&4\end{array}\right] ,
\]
we can specify the unique singular $K3$ surface.
We denote the singular $K3$ surface of discriminant $7$ by $X$ in this paper.

The $K3$ surface $X$ is just the same as the $K3$ surface $B/\Gamma_7^{\prime}$ described in Section 3.
We shall identify these $K3$ surfaces via the elliptic fibration $f_1$. 
In Figure \ref{Cox-graph}, three heptagons with vertices
\begin{gather*}
\begin{array}{ccccccc}
1&28&5&19&22&2&26,\\
14&10&24&16&12&15&9,\\
4&21&11&17&7&18&23
\end{array}
\end{gather*}
correspond to the reducible fibers, and the seven others correspond to the sections. 
Therefore the N{\'e}ron-Severi lattice $S_{B/\Gamma_7^{\prime}}$ is the overlattice of $U\oplus A_6^{\oplus 3}$ with a glue vector $(1,2,4)$ in the discriminant group $({\mathbb Z}/7{\mathbb Z})^3$ of $U\oplus A_6^{\oplus 3}$. 
By the Shioda-Tate formulae \cite{EMS}, this lattice is of rank $20$ and of discriminant $7$. Hence $B/\Gamma_7^{\prime}$ is isomorphic to $X$.

The orthogonal group ${\rm O}(T_X)$ of the transcendental lattice $T_X$ is isomorphic to $({\mathbb Z}/2{\mathbb Z})^2$. 
(This can be seen from the fact that $T_X$ contains only two $2$-roots.) 
The discriminant group $A_{T_X} \cong A_{S_X}$ of the transcendental lattice $T_X$ is isomorphic to $\mathbb{Z}/7\mathbb{Z}$ and the discriminant quadratic form is $\langle {4/7}\rangle$ on $A_{T_X}$. 
Especially, the orthogonal group ${\rm O}(q_{T_X})$ is $\{\pm 1\}$. 
There exists an automorphism of $X$ which acts as $-1$ on the ray ${\mathbb C}\omega_X$ of a holomorphic $2$-form $\omega_X$. 
For example, the deck transformation of a double plane ramified along the sextic curve in Section \ref{Naruki-results} (i.e. the inversion involution of the elliptic fibration of type $3I_7+3I_1$) acts on ${\mathbb C}\omega_X$ non-trivially since the quotient surface is rational. Especially, it acts as $-1$ on the discriminant group.
On the other hand, a cyclic subgroup of ${\rm O}(T_X){\cong}({\mathbb Z}/2{\mathbb Z})^2$ containing ${\mathbb Z}/2{\mathbb Z}$ generated by the deck transformation must be isomorphic to ${\mathbb Z}/2{\mathbb Z}$. 
Hence $H_X{\cong}{\mathbb Z}/2{\mathbb Z}$, and ${\rm Aut}(X)$ coincides with ${\rm Aut}(D(S_X))$. 

\vspace{2mm}
Instead of analyzing ${\rm Aut}(D(S_X))$ directly, we embed $S_X$ to $I{\!}I_{1,25}$. 
Then, by Nikulin's Theorem \ref{Nikulin}, the isometry of $S_X$ extends to that of $I{\!}I_{1,25}$. 
From Conway's Theorem \ref{Conway}, the desired group can be determined. 
The last step is described in the final Section. 

Consider the following vectors in $\Lambda$:
\begin{gather*}
Y=4\nu_0+\nu_{\Omega},\quad
Z=0,\quad
X=4\nu_{\infty}+\nu_{\Omega},\\
P=2\nu_{K_0},\quad
Q=4\nu_{\infty}+4\nu_{0},\quad
T=\nu_{\Omega}-4\nu_1,
\end{gather*}
where $K_0$ is an octad (i.e. an element of $S(5,8,24)$) which contains $0$ and does not contain $\infty$, $1$.
We take $\{ 0$ $2$ $3$ $4$ $5$ $15$ $20$ $22\}$ as $K_0$ throughout this paper.
Their corresponding Leech roots in $I{\!}I_{1,25}$ are:
\begin{gather*}
y=(2,1;Y),\quad
z=(-1,1;Z),\quad
x=(2,1;X),\\
p=(1,1;P),\quad
q=(1,1;Q),\quad
t=(1,1;T).
\end{gather*}
Then these six Leech roots generate an $A_6$ root lattice.
An $A_6$ root sublattice in $I{\!}I_{1,25}$ is unique up to the orthogonal group \cite[Chapter 23]{CS}. 
The orthogonal complement of this embedded root lattice $A_6$ is the lattice of signature $(1,19)$ with discriminant $7$. 
We denote the projection of a vector $x$ in $I{\!}I_{1,25}$ to $A_6^{\perp}$ by $x^{\prime}$ and the projection to $A_6^{\vee}$ by $x^{\prime\prime}$. 

The Leech roots orthogonal to $A_6$ are those twenty-eight Leech roots $r_i$ ($1\leq i\leq 28$) corresponding to the Leech lattice vectors $2\nu_{K_i}$, where $K_i$'s are the following octads in $\{ \infty ,0,1,\ldots ,22\}$ satisfying the conditions $\infty , 0\in K_i\not\ni 1$, $\sharp (K_0\cap K_i)=4$. 
These octads can be found in Todd's table \cite{Todd} as follows: 

\begin{gather*}
\begin{array}{@{\!\!\!\!\!}cccccccc@{\,}cccccccc}
K_{1\ }=\{\infty&0&  2&  3&  4&  8&  9& 21\} ,&
K_{2\ }=\{\infty&0&  2&  3&  6& 12& 16& 20\} , \\
K_{3\ }=\{\infty&0&  2&  3&  7& 11& 13& 15\} ,&
K_{4\ }=\{\infty&0&  2&  3& 10& 18& 19& 22\} , \\
K_{5\ }=\{\infty&0&  2&  4&  5&  6& 10& 11\} ,&
K_{6\ }=\{\infty&0&  2&  4&  7& 17& 18& 20\} , \\
K_{7\ }=\{\infty&0&  2&  4& 12& 14& 15& 19\} ,&
K_{8\ }=\{\infty&0&  2&  5&  7&  9& 12& 22\} , \\
\end{array}\\
\end{gather*}
\begin{gather*}
\begin{array}{@{\!\!\!\!\!}cccccccc@{\,}cccccccc}
K_{9\ }=\{\infty&0&  2&  5&  8& 13& 19& 20\} ,&
K_{10}=\{\infty&0&  2&  5& 15& 16& 18& 21\} , \\
K_{11}=\{\infty&0&  2&  6&  8& 15& 17& 22\} ,&
K_{12}=\{\infty&0&  2& 11& 14& 20& 21& 22\} , \\
K_{13}=\{\infty&0&  3&  4&  5& 12& 13& 18\} ,&
K_{14}=\{\infty&0&  3&  4&  6&  7& 14& 22\} , \\
K_{15}=\{\infty&0&  3&  4& 10& 15& 16& 17\} ,&
K_{16}=\{\infty&0&  3&  5&  6&  9& 15& 19\} , \\
K_{17}=\{\infty&0&  3&  5&  7& 10& 20& 21\} ,&
K_{18}=\{\infty&0&  3&  5&  8& 11& 16& 22\} , \\
K_{19}=\{\infty&0&  3&  8& 14& 15& 18& 20\} ,&
K_{20}=\{\infty&0&  3&  9& 13& 17& 20& 22\} , \\
K_{21}=\{\infty&0&  4&  5&  9& 14& 16& 20\} ,&
K_{22}=\{\infty&0&  4&  5& 17& 19& 21& 22\} , \\
K_{23}=\{\infty&0&  4&  6& 13& 15& 20& 21\} ,&
K_{24}=\{\infty&0&  4&  8& 10& 12& 20& 22\} , \\
K_{25}=\{\infty&0&  4&  9& 11& 15& 18& 22\} ,&
K_{26}=\{\infty&0&  5& 10& 13& 14& 15& 22\} , \\
K_{27}=\{\infty&0&  5& 11& 12& 15& 17& 20\} ,&
K_{28}=\{\infty&0&  7& 15& 16& 19& 20& 22\} . \\
\end{array}
\end{gather*}

Each pair of these twenty-eight Leech roots has the inner product $0$ or $1$.
Then Coxeter's graph \cite{Cox} appears as the graph whose vertices are these twenty-eight Leech roots and whose edges are the pair with the inner product $1$ (Figure \ref{Cox-graph}).

Conversely, twenty-eight Leech roots $r_i$'s span the overlattice of $U\oplus A_6^{\oplus 3}$ with a glue vector $(1,2,4)$ in the discriminant group $({\mathbb Z}/7{\mathbb Z})^3$. 
Since the discriminant $7$ contains no square factor, this lattice is primitive in $I{\!}I_{1,25}$. 
Hence this is $A_6^{\perp}$ in $I{\!}I_{1,25}$.

Thus $A_6^{\perp}$ in $I{\!}I_{1,25}$ can be identified with the N{\'e}ron-Severi lattice of the singular $K3$ surface of discriminant $7$. 
We denote by $C_i$ $(1\leq i\leq 28)$ the classes in $S_X$ corresponding to $r_i$. 
The twenty-eight $C_i$'s can be identified with the vertices of Coxeter's graph, i.e. the effective classes on $B/\Gamma_7^{\prime}$.

\begin{prop}\label{Weyl}
The projection $w^{\prime}$ of the Weyl vector $w=(1,0;0)$ to $(A_6^{\perp})^{\vee}$ is the sum of twenty-eight vertices $r_i$ of Coxeter's graph,that is, 
\[
w^{\prime}=\sum_{i=1}^{28}r_i.
\]
In particular, $w^{\prime}$ is in $A_6^{\perp}$, and $(w^{\prime})^2=28$.
\end{prop}
(Proof.) Since the twenty-eight Leech roots $r_i$ span $A_6^{\perp}$, it suffices to show $\langle w^{\prime}$, $r_j\rangle$ $=$ $\langle \sum r_i$, $r_j\rangle$ for any $1\leq{j}\leq{28}$.
The left-hand side is $\langle w^{\prime},r_j\rangle$ $=\langle w,r_j\rangle$ $=1$, this coincides with the right-hand side. 
The equality implies that $(w^{\prime})^2$ is twice the difference between the number of edges of Coxeter's graph and the number of vertices of Coxeter's graph. \hfill $\Box$
\begin{prop}\label{ample}
In the above identification of $A_6^{\perp}$ in $I{\!}I_{1,25}$ with $S_X$, $w^{\prime}$ is ample.
\end{prop}
(Proof.) Take an elliptic fibration of type $3I_7+3I_1$. 
Then Proposition \ref{Weyl} can be written as $w^{\prime}=3F+O+S_1+\cdots +S_6$ where $O$ is the zero section, $S_1,\ldots ,S_6$ are the non-zero sections, and $F$ is a general fiber. 
First, $(w^{\prime})^2=28>0$ holds.
To use Nakai's criterion of ampleness, we take an irreducible curve $C$.
If $\langle C,F\rangle =0$, then $C$ is a component of a fiber. 
Hence $\langle C,w^{\prime}\rangle$ $=1$ $>0$ holds. 
If $C$ is one of the sections $O$ or $S_i$, then $\langle C,w^{\prime}\rangle $ $=1$ $>0$ holds. 
In the other cases, $\langle C,O\rangle\geq 0$, $\langle C,S_i\rangle\geq 0$ and $\langle C,F\rangle >0$. 
Hence $\langle C,w^{\prime}\rangle >0$ holds. 
Therefore $w^{\prime}$ is ample.
\hfill $\Box$

\vspace{2mm}
We denote by $D^{\prime}$ the intersection of Conway's fundamental domain $D$ (Theorem \ref{Conway}) with $(A_6^{\perp})_{\mathbb R}$ in $(I{\!}I_{1,25})_{\mathbb R}$.
\begin{prop}{\rm (Borcherds \cite[Lemma 4.1]{Bor-Lor})}
The cone $D^{\prime}$ contains $w^{\prime}$. 
\end{prop}
(Proof.) 
To prove Proposition, it is enough to show that $w$ and $w^{\prime}$ are on a same side of the hyperplane $r^{\perp}$ for any Leech root $r$ $=$ $(-1-\lambda^2/2$, $1$; $\lambda )$ other than the simple roots $y$, $z$, $x$, $p$, $q$, $t$ of $A_6$. 
Since $\langle w^{\prime}$, $r\rangle$ $+\langle w^{\prime\prime}$, $r\rangle$ $=\langle w$, $r\rangle$ $=1$, it suffices to show $\langle w^{\prime\prime}$, $r\rangle\leq 0$. 
Since $\langle w^{\prime\prime}$, $y\rangle =$ $\langle w^{\prime\prime}$, $z\rangle =$ $\cdots =$ $\langle w^{\prime\prime}$, $t\rangle =$ $1$, the projection $w^{\prime\prime}$ is a linear combination of $y$, $z$, $x$, $p$, $q$, $t$ with negative coefficients, i.e. 
\[
w^{\prime\prime}=-3y-5z-6x-6p-5q-3t
\]
holds.
Then $\langle w^{\prime\prime}$, $r\rangle$ $\leq 0$ follows from the formula 
\[
\langle r_1,r_2\rangle =-2-\frac{(\lambda_1-\lambda_2)^2}{2}
\]
for two Leech roots $r_i=(-1-\lambda_i^2/2,1;\lambda_i)$ ($i=1,2$). 
\hfill $\Box$

Hence $D^{\prime}$ is a subcone of the ample cone $D(S_X)$ of $X$ under the above identification.

\begin{prop}
The stabilizer group of an $A_6$ diagram in ${\rm Aut}(D)$ is ${\rm PGL}_2(7)$.
\end{prop}
(Proof.)
The point-wise stabilizer of an $A_6$ lattice in ${\rm Aut}(D)$ is ${\rm PSL}_2(7)$ (see \cite[Chapter 10]{CS}). 
There exists a diagram reversing involution of an $A_6$ diagram (i.e. the opposition involution). 
Hence the stabilizer group as a set contains ${\rm PSL}_2(7)$ as a subgroup of index $2$.
On the other hand, the symmetry group of Coxeter's graph is ${\rm PGL}_2(7)$ (see \cite{Cox}), which contains ${\rm PSL}_2(7)$ as a subgroup of index $2$. 
\hfill $\Box$

\vspace{2mm}
\begin{lem}{\rm (Borcherds \cite[Lemma 4.2]{Bor-Lor})}
The faces of $D^{\prime}$ are the hyperplanes $r^{\prime\bot}$ perpendicular to $r^{\prime}$, where $r$ are the Leech roots such that $A_6\cup r$ generate a root lattice, and $r^{\prime}$ are their orthogonal projection to $S_X^{\vee}$. 
\end{lem}
(Proof.)
If the lattice generated by $A_6$ and $r$ is not negative definite, then $r^{\prime\bot}$$\cap$$P(S)$$=\emptyset$.
\hfill $\Box$

\begin{prop}\label{add-root}
Let $r$ be a Leech root that forms a root lattice $R$ together with $A_6$.
Then the following cases occur as the type of $R$:
\[
A_6\oplus A_1,\, A_7,\, D_7,\, E_7.
\]
The numbers of the Leech roots of the four types are $28$, $14$, $28$ and $56$, respectively.
\end{prop}
(Proof.)
A root lattice of rank $7$ which contains $A_6$ is one of the above four types. 
A Leech root generating one of the four types can be found in the Leech lattice as follows. 
In what follows, a capital $K$ denotes an octad and a small $k$ denotes an element of $\Omega$.

\begin{description}
\item[\fbox{$A_6{\oplus}A_1$}] A Leech root generating $A_6\oplus A_1$ together with $A_6$ is orthogonal to $A_6$. 
These Leech roots correspond to the Leech lattice elements $2\nu_{K^{(1)}}$, where $K^{(1)}$ is an octad such that $\infty ,0\in K^{(1)}\not\ni 1$ and $\sharp (K^{(1)}\cap K_0)=4$. The number of these Leech roots are 28. These Leech roots are the same as those illustrated in Figure \ref{Cox-graph}.
\item[\fbox{$A_7$ ($t$)}] A Leech root generating $A_7$ together with $A_6$ and an edge joining it with $t$ corresponds to the Leech lattice vector $2\nu_{K^{(2)}}$ where $K^{(2)}$ is an octad such that $\infty ,0,1\in K^{(2)}$ and $\sharp (K^{(2)}\cap K_0)=4$. The number of these Leech roots are 7.
\item[\fbox{$A_7$ ($y$)}] A Leech root generating $A_7$ together with $A_6$ and an edge joining it with $y$ corresponds to the Leech lattice vector $\nu_{\Omega}-2\nu_{K^{(3)}}+4\nu_{\infty}$ where $K^{(3)}$ is an octad such that  $\infty ,1\in K^{(3)}\not\ni 0$ and $K^{(3)}\cap K_0=\emptyset$. The number of these Leech root are $7$. 
\item[\fbox{$D_7$ ($q$)}] A Leech root generating $D_7$ together with $A_6$ and an edge joining it with $q$ corresponds to the Leech lattice vector $\nu_{\Omega}-4\nu_{k^{(1)}}$ where $k^{(1)}\in\Omega\setminus K_0$ and $k^{(1)}\not=\infty ,1$. The number of these Leech roots are $14$. 
\item[\fbox{$D_7$ ($z$)}] A Leech root generating $D_7$ together with $A_6$ and an edge joining it with $z$ corresponds to the Leech lattice vector $\nu_{\Omega} -2\nu_{K^{(4)}} +4\nu_{\infty} +4\nu_0 +4\nu_{k^{(2)}}$ where $k^{(2)}\in\Omega$ and $k^{(2)}\not=\infty ,0,1$ and $K^{(4)}$ is an octad such that $\infty, 0,1,k^{(2)}\in K^{(4)}$ and $K^{(4)}\cap K_0=\{ 0,k^{(2)}\}$. The number of these Leech roots are $14$. 
\item[\fbox{$E_7$ ($p$)}] A Leech root generating $E_7$ together with $A_6$ and an edge joining it with $p$ corresponds to the Leech lattice vector $2\nu_{K^{(5)}}$ where $K^{(5)}$ is an octad such that $\infty ,0\in K^{(5)}\not\ni 1$ and $\sharp (K^{(5)}\cap K_0)=2$. The number of these Leech roots are $28$. 
\item[\fbox{$E_7$ ($x$)}] A Leech root generating $E_7$ together with $A_6$ and an edge joining it with $x$ corresponds to the Leech lattice vector $\nu_{\Omega} -2\nu_{K^{(6)}} +4\nu_{\infty}$ where $K^{(6)}$ is an octad such that $0,1\in K^{(6)}\not\ni\infty$ and $\sharp (K^{(6)}\cap K_0)=2$. The number of these Leech roots are $28$. 
\end{description}

These numbers can be calculated in the same way as the case of the vertices of Coxeter's graph.
\hfill $\Box$
\vspace{2mm}

We call a Leech root generating $A_7$(resp. $D_7$, $E_7$)-type root lattice together with the $A_6$ root lattice an $A_7$(resp. $D_7$, $E_7$) Leech root, and we denote these Leech roots by $r_{A_7}$(resp. $r_{D_7}$, $r_{E_7}$). 
These Leech roots can be described in terms of the twenty-eight nodal curves $C_i$ $(1\leq i\leq 28)$ as follows. 

Each $A_7$ Leech root is orthogonal to twenty-four of $C_i$'s. 
It intersects with the four remaining $C_i$'s with multiplicity $1$. 
If we take $\{\infty$ $0$ $1$ $2$ $3$ $5$ $14$ $17\}$ as $K^{(2)}$, exceptional four nodal curves are $23$, $24$, $25$ and $28$ in the numbering of Figure \ref{Cox-graph}, and the black circles in Figure \ref{A_7}.
Each pair of these four vertices is of edge distance $4$. 
(Here, the {\it edge distance} means the minimal number of edges necessary for joining the two vertices.) 
For every vertex of Coxeter's graph, there exist two such $4$-tuples of vertices which contain that vertex.
These two $4$-tuples correspond to the $A_7$ Leech root joined to $A_6$ Dynkin diagram in opposite sides. 
The projections of $A_7$ Leech roots to $S_X$ are of norm $-8/7$.

Each $D_7$ Leech root is orthogonal to twenty-two of $C_i$'s. 
It intersects with the six remaining $C_i$'s with multiplicity $1$. 
If we take $6$ for $k^{(1)}$, the exceptional six nodal curves are $2$, $5$, $11$, $14$, $16$ and $23$ in Figure \ref{Cox-graph}, and the black circles in Figure \ref{D_7}. 
These six vertices in Coxeter's graph can be divided into three pairs of which are of edge distance $4$ and the remaining twelve pairs of which are of edge distance $3$. 
The projections of $D_7$ Leech roots to $S_X$ are of norm $-4/7$.

\begin{table}[bt]
\caption{The faces of $D^{\prime}$.}
\begin{center}
\begin{tabular}{c|cccc}
type $R$ & $A_6\oplus A_1$ & $A_7$ & $D_7$ & $E_7$ \\ \hline
$\sharp\{ r_R\}$ & $28$ & $14$ & $28$ & $56$ \\ 
$(r_R^{\prime})^2$ & $-2$ & $-8/7$ & $-4/7$ & $-2/7$ \\ 
$\langle w^{\prime},r_R^{\prime}\rangle$ & 1 & 4 & 6 & 7
\end{tabular}\\
\end{center}
\label{face-tab}
\end{table}

Each $E_7$ Leech root is orthogonal to twenty-one of $C_i$'s. It intersects with the seven remaining $C_i$'s with multiplicity $1$ (for example, $5$, $6$, $8$, $11$, $24$, $25$, $27$ in Figure \ref{Cox-graph} if we take $\{\infty$ $0$ $3$ $13$ $14$ $16$ $19$ $21\}$ as $K^{(5)}$ and the black circles in Figure \ref{E_7}). 
The projections of $E_7$ Leech roots are of norm $-2/7$.

Each $E_7$ Leech root defines the reflection on $S_X$. 
On the other hand, $D_7$ or $A_7$ Leech roots do not define the reflections on $S_X$. 
The symbol $s_{r_{D_7}^{\prime}}$ (resp. $s_{r_{A_7}^{\prime}}$) denotes the ${\mathbb Q}$-linear automorphisms of $(S_X)_{\mathbb Q}$ defined by the same rule as the reflection with respect to a $D_7$ (resp. $A_7$) Leech root. 

\begin{prop}\label{transitive}
The action of ${\rm PGL}_2(7)$ is transitive on the set of $A_7$ (resp. $D_7$ and $E_7$) Leech roots. 
\end{prop}
(Proof.)
This fact follows from the uniqueness result of the deep holes in the Leech lattice \cite[Chapter 23]{CS}, \cite[Lemma 6.1]{Bor-Lor}. \hfill $\Box$

These results are summarized in Table \ref{face-tab}.


\section{Some automorphisms of the singular $K3$ surface of discriminant $7$}

In this Section, we present several automorphisms of the singular $K3$ surface of discriminant $7$ associated with some elliptic fibrations.

\vspace{2mm}
For each $E_7$ (resp. $D_7$, $A_7$)-type face of the polyhedral cone $D^{\prime}$ in Proposition \ref{add-root}, there is an elliptic fibration of type $6I_1+I_{18}$ (resp. $4I_1+I_{12}+I\!V^*$, $2I_1+I_8+2I_1^*$) in Kodaira's notation. 
Weierstrass equations of these elliptic fibrations were obtained by Harrache and Lecacheux \cite{HL}.

\begin{figure}[tb]
\caption{$I_{18}$-fiber in Coxeter's graph.}
\begin{center}
\includegraphics[height=240pt]{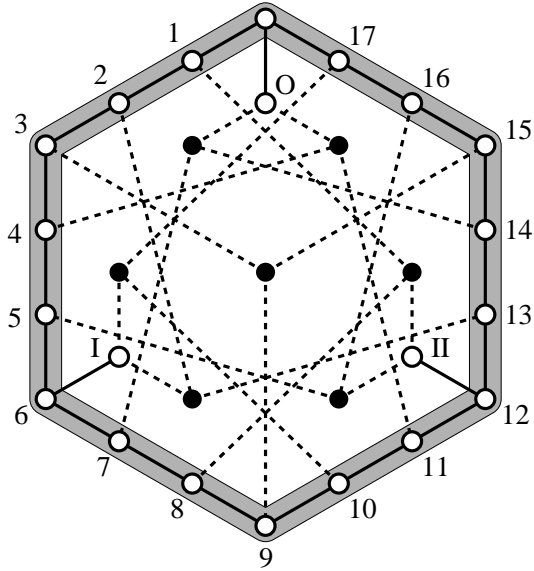}\\
\end{center}
\label{E_7}
\end{figure}

\begin{figure}[tb]
\caption{($I_{12}+I\!V^*$)-fibers in Coxeter's graph.}
\begin{center}
\includegraphics[height=240pt]{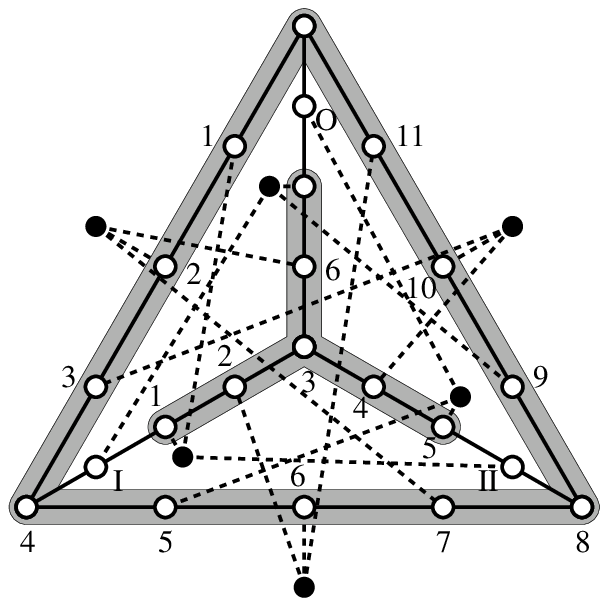}\\
\end{center}
\label{D_7}
\end{figure}

\begin{figure}[tb]
\caption{($I_8+2I_1^*$)-fibers in Coxeter's graph.}
\begin{center}
\includegraphics[height=240pt]{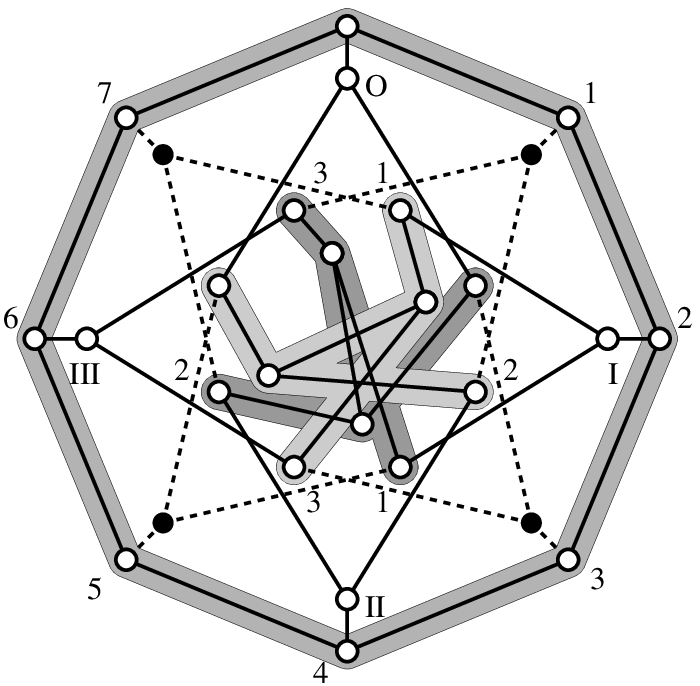}\\
\end{center}
\label{A_7}
\end{figure}

These fibrations can be found in Coxeter's graph as given in Figures \ref{E_7}, \ref{D_7} and \ref{A_7}. 
If one removes the curves intersecting an $E_7$-root with multiplicity $1$ (these are given at the end of Section 4, and denoted by black circles in Figure \ref{E_7}), a Dynkin diagram of type $A_{17}^{\sim}$ can be found in the remaining vertices. 
If one removes the curves intersecting a $D_7$-root with multiplicity $1$ (black circles in Figure \ref{D_7}), there remains the graph analogous to the skeleton of a tetrahedron. 
Then four Dynkin diagrams of type $E_6^{\sim}$ can be found at the four vertices of the tetrahedron. 
If a Dynkin diagram of type $E_6^{\sim}$ is fixed, then a Dynkin diagram of type $A_{11}^{\sim}$ can be found in the remains.
Similarly, a Dynkin diagram of type $A_7^{\sim}{\oplus}D_5^{\sim\oplus 2}$ can be found in the case of an $A_7$ root in three ways. 
One of these Dynkin diagrams are indicated by gray shadows in Figures. 
(In Figures, Roman numerals indicate the sections, whereas Arabic numerals give the numbering of the basis of these root lattices, which differ from the numbering in Figure \ref{Cox-graph}.) 
We fix one elliptic fibration $f_R$ to each $R$ Leech root. 

These fibrations have the trivial lattices of rank $19$, and their essential lattices are of rank $1$. 
The projection $r^{\prime}_R$ ($R=E_7$, $D_7$ or $A_7$) of an $R$ Leech root $r_R$ to $S_X^{\vee}$ is orthogonal to the trivial lattice of the corresponding fibration $f_R$, therefore $r_R^{\prime}$ generates the essential lattice ${\rm Ess}_{\mathbb Q}$ over ${\mathbb Q}$. 
We denote the inversion involutions of these three fibrations by $\iota_{E_7}$, $\iota_{D_7}$, or $\iota_{A_7}$ according to their types. 
Then, $\iota_R (r_R^{\prime})=-r_R^{\prime}$ holds by Proposition \ref{inv} ($R=E_7$, $D_7$ or $A_7$).


As shown in Figure \ref{E_7}, $E_7$ Leech roots correspond to the $I_{18}$-fibrations and $3$-torsion sections.
For each $r_{E_7}$, seven vectors out of twenty-eight vertices of Coxeter's graph have inner product $1$, and the other vectors are orthogonal to $r_{E_7}$. 
Therefore the stabilizer group of an $E_7$ face in ${\rm PGL}_2(7)$ is isomorphic to the dihedral group $\mathfrak{D}_3$ of degree $3$.
The stabilizer group of a $D_7$-face in ${\rm PGL}_2(7)$ is isomorphic to the alternating group $\mathfrak{A}_4$.
The stabilizer group of an $A_7$-face is isomorphic to the symmetric group $\mathfrak{S}_4$. 


\begin{prop}\label{q-refl}
The inversion involution $\iota_R$ ($R=E_7$, $D_7$ or $A_7$) acts as the reflection with respect to $r_R^{\prime}$ on the plane spanned by $w^{\prime}$ and $r_R^{\prime}$.\\
That is, 
\begin{eqnarray*}
\iota_R(r_R^{\prime})&=&s_{r_R^{\prime}}(r_R^{\prime})=-r_R^{\prime}, \\
\iota_R(w^{\prime})&=&s_{r_R^{\prime}}(w^{\prime})
\end{eqnarray*}
where $R=E_7$, $D_7$, or $A_7$.
\end{prop}
(Proof.)
The first row was already shown in Proposition \ref{inv}.
The projection $w^{\prime}$ of Weyl vector $w$ to $S_X$ is decomposed as
\[
w^{\prime}=w^{\prime}_{\rm triv}+w^{\prime}_{\rm ess}
\]
according to the decomposition $S_{\mathbb Q} \cong {\rm Triv}_{\mathbb Q}\oplus{\rm Ess}_{\mathbb Q}$ for the elliptic fibration $f=f_R$.

By Proposition \ref{inv}, $\iota_R(w^{\prime}_{\rm ess}) = s_{r_R^{\prime}}(w^{\prime}_{\rm ess})$ holds. 
Hence there is nothing to prove on the essential part. 
On the other hand, $s_{r_R}(w^{\prime}_{\rm triv})=w^{\prime}_{\rm triv}$ holds. 
Hence it suffices to prove that $\iota_R(w^{\prime}_{\rm triv}) = w^{\prime}_{\rm triv}$. 
We can verify this equality by projecting $w^{\prime}_{\rm triv}$ to each component of the decomposition $U_f\oplus \left( \bigoplus \Theta_{i}\right)$ of ${\rm Triv}$. 

\vspace{2mm}
For example, we see the case of a $D_7$ Leech root.
In this case, the trivial lattice is $U_f\oplus A_{11}\oplus E_6$. 
We check the action of $\iota_{D_7}$ on each component. 


The classes of the zero section and a general fiber are fixed by the inversion. 
Hence the component $U_f$ is fixed. 
To consider the projection of $w^{\prime}$ to the component $A_{11}$, we decompose the summation $w^{\prime}=\sum_{i=1}^{28} C_i$ (Proposition \ref{Weyl}) as 
\[
w^{\prime}=\sum_{A_{11}^{\sim}}C_i+\sum_{E_6^{\sim}}C_i+\sum_{\mbox{sections}}C_i
  +\sum_{\mbox{black circles}}C_i.
\]
The sum $\sum_{A_{11}^{\sim}}C_i$ of the twelve roots of components in the singular fiber of type $A_{11}^{\sim}$ is in $U_f$.
Hence its projection to $A_{11}$ is zero. 
The $(-2)$-roots $C_i$ of the classes of components in the singular fiber of type $E_6^{\sim}$ are orthogonal to $A_{11}$, and the projection of $\sum_{E_6^{\sim}}C_i$ to $A_{11}$ is zero. 
The class of the zero section is orthogonal to $A_{11}$. 
We use the numbering in Figure \ref{D_7} for the basis $e_i$ of $A_{11}$, and let ${\rm Gram}_{A_{11}}$ be the Gram matrix of $A_{11}$ with respect to this basis. 
The projections of the two roots of non-zero sections are expressed as ${\rm Gram}_{A_{11}}^{-1}(e_4)$ and ${\rm Gram}_{A_{11}}^{-1}(e_8)$.
Then the sum ${\rm Gram}_{A_{11}}^{-1}(e_4)$$+{\rm Gram}_{A_{11}}^{-1}(e_8)$ is fixed by the inversion. 
The projections of the six roots intersecting $r_{D_7}$ can be expressed as
\begin{gather*}
{\rm Gram}_{A_{11}}^{-1}e_1,\quad
{\rm Gram}_{A_{11}}^{-1}(e_2+e_7),\quad
{\rm Gram}_{A_{11}}^{-1}(e_3+e_{10}),\\
{\rm Gram}_{A_{11}}^{-1}e_5,\quad
{\rm Gram}_{A_{11}}^{-1}(e_6+e_{11}),\quad
{\rm Gram}_{A_{11}}^{-1}e_9.
\end{gather*}
The sum of these six vectors is
\[
{\rm Gram}_{A_{11}}^{-1}(e_1+e_2+e_3+e_5+e_6+e_7+e_9+e_{10}+e_{11}).
\]
This is fixed by the inversion. (The six roots not orthogonal to $r_{D_7}$ are joined by edge to the singular fiber of type $A_{11}^{\sim}$ not symmetrically for inversion. But their sum is symmetric.)

The component $E_6$ can be treated like as the component $A_{11}$. 
The sum of the roots orthogonal to the $D_7$ Leech root is fixed by the inversion. 
Under the numbering in Figure \ref{D_7} for the basis $e_i$ of $E_6$, the projections of the remaining six roots to the $E_6$ component are 
\begin{gather*}
{\rm Gram}_{E_6}^{-1}e_1,\,
{\rm Gram}_{E_6}^{-1}e_2,\,
{\rm Gram}_{E_6}^{-1}e_3,\,
{\rm Gram}_{E_6}^{-1}e_4,\,
{\rm Gram}_{E_6}^{-1}e_5,\,
{\rm Gram}_{E_6}^{-1}e_6 .
\end{gather*}
Therefore the sum of these six vectors is fixed by the inversion. 

The cases of an $E_7$ and an $A_7$ Leech root can be verified similarly.
\hfill $\Box$


\section{A system of generators of the automorphism group}

To prove that the automorphism group of $X$ is generated by the inversion involutions constructed in previous section together with ${\rm PGL}_2(7)$, the same argument as Kond{\=o} \cite{Kondo} can be applied. 

\begin{thm}\label{gen}
Let $X$ be the singular $K3$ surface of discriminant $7$. 
And let $N$ be the subgroup of ${\rm Aut}(X)$ generated by the inversion involutions $\iota_R$ ($R=E_7$, $D_7$, and $A_7$) defined in the previous section.
Then the automorphism group of $X$ is a split extension of the finite group ${\rm PGL}_2(7)$ by $N$.
\end{thm}
(Proof.)
Let $\varphi$ be an element of ${\rm Aut} (X)$. Then it suffices to show that there exists an automorphism $\psi$ in $N$ such that $\psi\circ\varphi\in {\rm Aut} (D^{\prime})$.
Let $w^{\prime}$ be the projection of the Weyl vector $w=(1,0;0)$ to $S_X$. Take an element $\psi_1$ of $N$ that minimizes $\langle \psi_1\circ\varphi (w^{\prime}),w^{\prime}\rangle$.
Hence, for the projection $r_R^{\prime}$ ($R=E_7$, $D_7$ or $A_7$) to $(S_X)_{\mathbb Q}$ of any $R$ Leech root $r_R$ of Proposition \ref{add-root}, 
\begin{eqnarray*}
\langle\psi_1\circ\varphi (w^{\prime}),w^{\prime}\rangle
&\leq &
\langle \iota_R\circ\psi_1\circ\varphi (w^{\prime}),w^{\prime}\rangle
=
\langle\psi_1\circ\varphi (w^{\prime}),\iota_R(w^{\prime})\rangle \\
&=&
\langle\psi_1\circ\varphi (w^{\prime}),w^{\prime}\rangle -2
\frac{\langle w^{\prime},r_R^{\prime}\rangle}{\langle r_R^{\prime},r_R^{\prime}\rangle}
\langle\psi_1\circ\varphi (w^{\prime}),r_R^{\prime}\rangle
\end{eqnarray*}
holds, because the automorphism $\iota_R$ acts as the reflections to $w^{\prime}$. Here 
\[
\langle w^{\prime},r_R^{\prime}\rangle > 0,\quad
\langle r_R^{\prime},r_R^{\prime}\rangle < 0
\]
hold. 
Therefore we obtain the inequalities $\langle\psi_1\circ\varphi (w^{\prime}),r_R^{\prime}\rangle\geq 0$ for all $r_R$. 
Hence $\psi_1\circ\varphi (w^{\prime})$ is contained in $D^{\prime}$.
\hfill $\Box$

{\sc Graduate School of Mathematics, Nagoya University, Chikusa-Ku Nagoya 464-8602 Japan}\\
{\it E-mail address}: {\tt d07001n@math.nagoya-u.ac.jp}

\end{document}